\documentclass[12pt]{amsart}
\usepackage{amssymb}
\begin{document}

\def\dbl{[\hskip -1pt[}
\def\dbr{]\hskip -1pt]}
\newcommand{ \al}{\alpha}
\newcommand{\ab}[1]{\vert z\vert^{#1}}
\title
{Higher order
symmetries
of
real
hypersurfaces in $\mathbb C^3$}

\author{ martin kolar and francine meylan }

\address{M. Kolar: Department of Mathematics and Statistics, Masaryk University,
Kotlarska 2,  611 ~37
 Brno, Czech Republic}

\email{mkolar@math.muni.cz}
\thanks{The first author was supported by the project CZ.1.07/2.3.00/20.0003
of the Operational Programme Education for Competitiveness of the Ministry
of Education, Youth and Sports of the Czech Republic.}

\address{ F. Meylan: Department of Mathematics,
University of Fribourg, CH 1700 Perolles, Fribourg}

\email{francine.meylan@unifr.ch}
\thanks{The second author was supported by Swiss NSF Grant 2100-063464.00/1 }

\begin{abstract} 
We study nonlinear automorphisms
of Levi degenerate hypersurfaces of finite multitype.
 By results of \cite{KMZ}, the Lie algebra of infinitesimal CR automorphisms 
may contain 
a graded component 
 consisting of nonlinear vector fields of arbitrarily high degree,
which has no analog 
in the classical Levi nondegenerate case, or in the case of finite type hypersurfaces in $\mathbb C^2$. 
We analyze this phenomenon for hypersurfaces of finite Catlin multitype in complex dimension three.
The results provide a complete classification of such manifolds.
As a consequence, we
show on which hypersurfaces  
2-jets
are not sufficient to determine an automorphism.
%
The results also confirm a conjecture about the origin of nonlinear automorphisms of Levi degenerate hypersurfaces,
 formulated by the first author (AIM 2010).
\end{abstract} 


\def\Label#1{\label{#1}}
\def\1#1{\ov{#1}}
\def\2#1{\widetilde{#1}}
\def\6#1{\mathcal{#1}}
\def\4#1{\mathbb{#1}}
\def\3#1{\widehat{#1}}
\def\K{{\4K}}
\def\LL{{\4L}}
\def\H{{\4H}}
\def\C{{\4C}}
\def\R{{\4R}}
\def \MM{{\4M}}
\def\La{\Lambda}
\def\la{\lambda}
\def\Re{{\sf Re}\,}
\def\Im{{\sf Im}\,}

\numberwithin{equation}{section}
\def\s{s}
\def\k{\kappa}
\def\ov{\overline}
\def\span{\text{\rm span}}
\def\ad{\text{\rm ad }}
\def\tr{\text{\rm tr}}
\def\xo {{x_0}}
\def\Rk{\text{\rm Rk\,}}
\def\sg{\sigma}
\def \emxy{E_{(M,M')}(X,Y)}
\def \semxy{\scrE_{(M,M')}(X,Y)}
\def \jkxy {J^k(X,Y)}
\def \gkxy {G^k(X,Y)}
\def \exy {E(X,Y)}
\def \sexy{\sccenikrE(X,Y)}
\def \hn {holomorphically nondegenerate}
\def\hyp{hypersurface}
\def\prt#1{{\partial \over\partial #1}}
\def\det{{\text{\rm det}}}
\def\wob{{w\over B(z)}}
\def\co{\chi_1}
\def\po{p_0}
\def\fb {\bar f}
\def\gb {\bar g}
\def\Fb {\ov F}
\def\Gb {\ov G}
\def\Hb {\ov H}
\def\zb {\bar z}
\def\wb {\bar w}
\def \qb {\bar Q}
\def \t {\tau}
\def\z{\chi}
\def\w{\tau}
\def\Z{\zeta}
\def\phi{\varphi}
\def\eps{\varepsilon}

\def \T {\theta}
\def \Th {\Theta}
\def \L {\Lambda}
\def\b {\beta}
\def\a {\alpha}
\def\o {\omegacenik}
\def\l {\lambda}

\def \im{\text{\rm Im }}
\def \re{\text{\rm Re }}
\def \Char{\text{\rm Char }}
\def \supp{\text{\rm supp }}
\def \codim{\text{\rm codim }}
\def \Ht{\text{\rm ht }}
\def \Dt{\text{\rm dt }}
\def \hO{\widehat{\mathcal O}}
\def \cl{\text{\rm cl }}
\def \bR{\mathbb R}
\def \bS{\mathbb S}
\def \bK{\mathbb K}
\def \bD{\mathbb D}
\def \bC{\mathbb C}
\def \C{\mathbb C}
\def \N{\mathbb N}
\def \bL{\mathbb L}
\def \bZ{\mathbb Z}
\def \bN{\mathbb N}
\def \scrF{\mathcal F}
\def \scrK{\mathcal K}
\def \mc #1 {\mathcal {#1}}
\def \scrM{\mathcal M}
\def \cR{\mathcal R}
\def \scrJ{\mathcal J}
\def \scrA{\mathcal A}
\def \scrO{\mathcal O}
\def \scrV{\mathcal V}
\def \scrL{\mathcal L}
\def \scrE{\mathcal E}
\def \hol{\text{\rm hol}}
\def \aut{\text{\rm aut}}
\def \Aut{\text{\rm Aut}}
\def \J{\text{\rm Jac}}
\def\jet#1#2{J^{#1}_{#2}}
\def\gp#1{G^{#1}}
\def\gpo{\gp {2k_0}_0}
\def\emmp {\scrF(M,p;M',p')}
\def\rk{\text{\rm rk\,}}
\def\Orb{\text{\rm Orb\,}}
\def\Exp{\text{\rm Exp\,}}
\def\Span{\text{\rm span\,}}
\def\d{\partial}
\def\D{\3J} 
\def\pr{{\rm pr}}

\def \CZZ {\C \dbl Z,\zeta \dbr}
\def \D{\text{\rm Der}\,}
\def \Rk{\text{\rm Rk}\,}
\def \CR{\text{\rm CR}}
\def \ima{\text{\rm im}\,}
\def \I {\mathcal I}

\def \M {\mathcal M}
\def\5#1{\mathfrak{#1}}

\newtheorem{Thm}{Theorem}[section] 
\newtheorem{Cor}[Thm]{Corollary}
\newtheorem{Pro}[Thm]{Proposition}
\newtheorem{Lem}[Thm]{Lemma}
\newtheorem{Conj}[Thm]{Conjecture}

\theoremstyle{definition}\newtheorem{Def}[Thm]{Definition}

\theoremstyle{remark}
\newtheorem{Rem}[Thm]{Remark}
\newtheorem{Exa}[Thm]{Example}
\newtheorem{Exs}[Thm]{Examples}

\def\bl{\begin{Lem}}
\def\el{\end{Lem}}
\def\bp{\begin{Pro}}
\def\ep{\end{Pro}}
\def\bt{\begin{Thm}} 
\def\et{\end{Thm}}
\def\bc{\begin{Cor}}
\def\ec{\end{Cor}}
\def\bd{\begin{Def}}
\def\ed{\end{Def}}
\def\br{\begin{Rem}}
\def\er{\end{Rem}}
\def\be{\begin{Exa}}
\def\ee{\end{Exa}}
\def\bpf{\begin{proof}}
\def\epf{\end{proof}}
\def\ben{\begin{enumerate}}
\def\een{\end{enumerate}}

\newcommand{\zz}{(z,\bar z)}  

\maketitle 

\section{Introduction}


One of the central problems in CR geometry  is  the classification of real hypersurfaces in $\mathbb C^n$, 
up to  biholomorphic equivalence. A complete solution of this problem should also lead to
 a complete 
understanding of automorphism groups of such manifolds.

When the hypersurface is Levi nondegenerate, the problem is 
well
understood thanks  to   the classical work of  Chern and Moser
\cite{CM}. In particular, 
the infinitesimal  CR  automorphisms of such manifolds form a graded Lie algebra
with at most 5 components. Moreover, by results of Kruzhilin and Loboda (\cite{KL}),  if a strongly
pseudoconvex  hypersurface 
is not equivalent to the sphere, there are at most 3 graded components, and 
all infinitesimal automorphisms are linear in appropriate coordinates.
For the sphere itself, the coefficients are at most quadratic, 
which implies  2-jet determination in general.

Similar results were obtained for hypersurfaces of finite type in $\mathbb C^2$. In particular, 
the same 2-jet determination result holds
(see  \cite{ELZ2}, \cite{KM3}).

In a  recent paper 
\cite{KMZ}, the same problem is considered for Levi degenerate hypersurfaces in $\mathbb C^n$
 with weighted homogeneous polynomial models, 
 which replace 
 the model hyperquadric from the 
nondegenerate case. 
The results
describe possible structures of infinitesimal CR automorphism  algebras
for hypersurfaces of finite Catlin multitype.
Compared to the Levi nondegenerate case, there are in general 6 possible 
components.
The new phenomenon is the existence
of nonlinear infinitesimal CR automorphisms in the complex tangential variables,
which are of arbitrarily high degree in general.

Our aim in this paper  is to analyze this phenomenon
and provide a complete description  of 
hypersurfaces of finite Catlin multitype  in $\mathbb C^3$ which admit such  automorphisms.

Let us recall that the Catlin multitype is an important CR invariant which 
Catlin introduced to prove subelliptic estimates on pseudoconvex domains (\cite{C},\cite{C1}).
The definition of multitype  was extended to the general case (not necessarily pseudoconvex)  in \cite{Ko1}.
%
It provides a natural setting for an extension of the Chern-Moser theory
to  degenerate manifolds
(\cite{KMZ}).
We consider 
 a weighted homogeneous 
model 
of finite Catlin multitype that is holomorphically nondegenerate.
Let 
\begin{equation}\Label{model}
M_{P}:=\{\Im w = P(z,\bar z)\}, \quad (z,w)\in \C^2\times\C,
\end{equation}
where $P$ is a real valued weighted homogeneous polynomial with respect to the multitype weights $\mu_1, \mu_2$ 
(see Section 2 for the needed definitions).



As proved in \cite{KMZ}, the   Lie algebra of infinitesimal CR automorphisms
$\5g=\aut(M_{P},0)$ of $M_{P}$ admits the  weighted decomposition given
by
\begin{equation}
\5g = \5g_{-1} \oplus \bigoplus_{j=1}^{2}\5g_{- \mu_j} \oplus \5g_{0}
\oplus \5g_c \oplus  \5g_{n}
 \oplus \5g_{1}.
 \label{muth}
\end{equation}
Here $\5g_{c}$ contains  vector fields commuting with $W$ and 
 $\5g_{n}$ contains vector fields not commuting with $W$, both of weight $\mu \in (0,1)$ (see \cite{KMZ} for 
 more details). Notice that $\5g_{-1}$  contains $W=\d_{w}$  and $\5g_{0}$ contains the weighted Euler field, hence they
are always nontrivial.
 A complete description of $\5g_1$ is also contained in \cite{KMZ}.  

\begin{Rem}\label{rustu} By the results of \cite{KMZ}, the elements of $\5g_{n}$ and $\5g_{1}$
are determined by ordinary 2-jets, hence {\it higher order} infinitesimal automorphisms may occur only  when   
$\5g_{c}$ is nontrivial.
\end{Rem}

In this paper, we  study  all hypersurfaces whose model has nontrivial  $\5g_{c}$. 
Our results confirm a conjecture about the origin of nonlinear automorphisms of Levi degenerate hypersurfaces
 formulated by the first author (see the 2010 AIM list of problems
http://www.aimath.org/WWN/crmappings/crmappings.pdf): 
$M_P$ has a nonlinear symmetry if and only if there is a holomorphic mapping $f$ of $M_p$ into a hyperquadric 
in $\mathbb C^K$ 
and a symmetry of the hyperquadric, which
is $f$-related to the automorphism of $M_P$. 

 Note that mappings of CR manifolds into hyperquadrics have been studied intensively in recent years (see e.g. 
\cite{BEH}, \cite{EHZ}).
Here we ask in addition 
compatibility with  a symmetry of the hyperquadric.
Let us remark that analysing $\5g_n$ requires completely different techniques, and is the subject of a forthcoming paper \cite{KMp}

In order to  formulate  our first result, let us 
recall that two vector fields $X_1$ and $X_2$  are  $f$-related (or compatible by $f$) if $f_*(X_1) = X_2.$
\bt
\label{3piva}
Let $M_{P}$ be a holomorphically nondegenerate hypersurface  given by \eqref{model}. Assume that $  \dim \5g_c > 0$
and let $Y \in \5g_c$ be a nonzero vector field.
Then there exist an integer  $K \ge 3$  and 
a holomorphic  polynomial mapping $f:\mathbb C^3 \longrightarrow \mathbb C^K$
 which maps  $M_{P}$ into a hyperquadric $H\subseteq \mathbb C^K,$  such that 
 $Y$ is $f$-related with an infinitesimal CR automorphism of $H$.
\et

The proof is based on an explicit complete description of models with nontrivial $\5g_c$, for which 
we need the following definition.
%
%
\bd
Let $Y$ be a weighted homogeneous vector field. 
A pair of finite sequences   of holomorphic  weighted homogeneous polynomials $\{U^1, \dots, U^n\}$ and  $\{V^1, \dots, V^n\}$ is called a symmetric pair of $Y-$chains if 
\begin{equation}
Y(U^n)=0, \ \ Y(U^j) =c_j U^{j+1}, \ \ j=1, \dots, n-1,
\label{yuv}
\end{equation}
\begin{equation}
Y(V^n)=0, \ \ Y(V^j) =d_j V^{j+1}, \ \ j=1, \dots, n-1,
\end{equation}
where $c_j$, $d_j$ are non zero complex constants, which satisfy
\begin{equation} c_j= - \bar d_{n-j}. \label{dcj}
\end{equation}
If the two sequences are identical 
we say that $\{U^1, \dots, U^n\}$ is a symmetric $Y$ - chain.
\ed

\be
Let  $$Y=i{z_2}^l\dfrac{\partial}{\partial z_1}.
  $$
Then the pair $\{U^1, U^2 \}= \{z_1, z_2^l\}$ is a symmetric $Y-$chain, since $Y(U^2) =0$ and $Y(U^1) = i U^2.$ Then for the hypersurface given by 
$$  \Im w = \Re U^1 \overline{ U^2}  = \Re z_1{\overline{z_2}}^l$$ we have  $Y \in \5g_c$.

\ee

The following result shows that in general the elements of $\5g_c$  arise in an analogous way.

\bt \label{pivo} Let 
 $M_P$  be  given by  \eqref{model} admitting a nontrivial 
 $Y\in \5g_c$. 
Then $P$ can be decomposed in the following way

\begin{equation}P=\sum_{j=1}^M T_j,\label{Barou19}
\end{equation} where each $T_j$ is given by 
\begin{equation}\label{Barou20}
 T_j =\Re (\sum_{k=1}^{N_j}
  {U_j^{k}}
  {\overline {{V_j^{N_j -k +1}}}}), 
  \end{equation}
where $\{{ {{U_j^{1}}, \dots, {U_j^{N_j}} }}\}$  and  $\{{ {{V_j^{1}}, \dots, {V_j^{N_j}} }}\}$ are a symmetric pair of $Y-$ chains.


Conversely, if $Y$ and $P$ satisfy  \eqref{yuv} --
\eqref{Barou20},
then $Y \in \5g_c$.
\et
\br
It is immediate to see that  $Y$ is uniquely and explicitely determined by $P$. More precisely, since $M_P$ is holomorphically nondegenerate, 
at least one of the $T_j$ has length $N_j \ge 2. $ For such a $T_j$ we have 
$$Y = \frac{c_{N_j-1} U_j^{N_j-1}}{\Delta} \left(-\frac{\partial U^{N_j}_j}{\partial z _2}\d_{z_1}   +  \frac{\partial U_j^{N_j}}{\partial z_1}\d_{z_2}\right), $$
where $\Delta $ is the Jacobian of $\{U_j^{N_{j-1}}, U_j^{N_j}\}$.  Hence for a given hypersurface 
the results also provide a simple constructive tool to determine $\5g_c$. Moreover, this also shows that the real dimension of $\5g_c$
is at most one. 
\er
Examples of symmetric chains of arbitrary length are described at the end of Section 3.
Using  Remark \ref{rustu}, we obtain
 
 \bt \label{7piv}
 Let $M$ be an arbitrary smooth hypersurface of finite Catlin multitype. If its model is holomorphically 
 nondegenerate and  not biholomorphically equivalent to one of 
 the form described in Theorem \ref{pivo},
 then 
 the automorphisms of $M$ are determined by their  2-jets.
 \et

%

The paper is organized as follows. Section 2 contains the necessary definitions used in the rest of the paper.
Section 3 deals with the $\5g_c$ component of the algebra $\aut(M_{P},0)$.
Section 4 completes the proofs of the above theorems.

\section{Preliminaries}

Let  $M \subseteq \mathbb C^{3}$ be a smooth hypersurface,
and $p \in M $ be a  point of { finite type} $m \ge 2$  in the
sense of Kohn and Bloom-Graham (\cite{BG},\cite{BG1},  \cite{K}).

\noindent We  consider 
local holomorphic coordinates $(z,w)$ vanishing at $p$,
where $z =(z_1, z_2)$ and  $z_j = x_j + iy_j$, $j=1,2$, and
$w=u+iv$. The hyperplane $\{ v=0 \}$ is assumed to be tangent to
$M$ at $p$, hence  $M$  is described near $p=0$ as the graph of a uniquely
determined real valued function
\begin{equation} v = \phi(z_1,z_2, \bar z_1,z_2,  u), \ d\phi(0) =0.
\label{vp}
\end{equation}

We can assume 
that (see e.g. \cite{BG})
\begin{equation}\label{fifi}
\phi(z_1,z_2,  \bar z_1,\bar z_2,  u)=P_m(z, \bar z) +o(u,|z|^m),
\end{equation}
where $P_m(z, \bar z)$ is a nonzero homogeneous polynomial of degree $m$ without pluriharmonic terms.

Recall that the definition of multitype involves   rational  weights associated  to the variables
$w, z_1, z_2$.
%
%
%
%
%
%
%
 The
variables $w$, $u$ and $v$ are given weight one, { reflecting} our choice
of  tangential and normal variables.
The complex tangential variables $(z_1, z_2)$  are treated according to
the following definitions (for more details, see \cite{Ko1}).

\bd
A weight is a pair of nonnegative
 rational numbers $\La = (\la_1,
\la_2)$, where $0 \leq\la_j\leq \frac12$, and $\la_1 \ge
\la_{2}$.
\ed

Let $\La = (\la_1,
\la_2)$ be a weight, and   $\al=(\al_1, \al_2), \ $
 $\ \beta=(\beta_1,\beta_2) $   be  multiindices.
The weighted degree $\kappa$ of a monomial $$q(z, \bar z,u)=c_{\al
\beta l}z^{\al}\bar z^\beta u^{l} , \ l \in \mathbb N,$$ is defined as
$$ \kappa:=
l +  \sum_{i=1}^2 (\al_i + \beta_i ) \la_i.$$

%

{A polynomial $Q(z, \bar z, u)$  is weighted homogeneous
of weighted degree $\kappa$ if it is a sum of
 monomials of weighted degree $\kappa$.}


For a weight  $\La$,  the weighted length of a multiindex $\al = (\al_1, \al_2)$ is
defined by

$$\vert \alpha\vert_{\La} := \la_1 \al_1 + \la_2 \al_2.$$

Similarly, if $\al = (\al_1, \al_2)$ and  $\hat \al =
(\hat \al_1, \hat \al_2)$ are two multiindices, the weighted
length of the  pair $(\al, \hat \al)$ is
$$\vert (\alpha,\hat \al) \vert_{\La} := \la_1 (\al_1 +\hat \al_1) +
\la_2 (\al_2 + \hat \al_2).$$


\bd{ A weight $\La$ will be called distinguished for $M$ if there exist
local holomorphic coordinates $(z,w)$ in which the defining equation of $M$ takes form
\begin{equation} v = P\zz + o_{\La}(1),
\label{1}
\end{equation}
where $P\zz$ is a nonzero $\La$ - homogeneous polynomial of
weighted degree $1$ without pluriharmonic terms, and $o_{\La}(1)$
denotes a smooth function whose derivatives of weighted order less than or equal to
one vanish.}
\ed

The fact that distinguished weights do exist follows from \eqref{fifi}.
For these coordinates $(z,w),$ we have
$$\Lambda=(\dfrac{1}{m}, \dfrac{1}{m}).$$

In the following we shall consider the standard lexicographic order on the set
of pairs.

We recall the following definition (see \cite{C}).

\bd \label{2.6}
Let  $\Lambda_M = (\mu_1, \mu_2)$  be the infimum of all possible
distinguished weights  $\Lambda$ with respect to the lexicographic order.
The multitype of $M$ at $p$ is defined to be the pair $$(m_1,
 m_2),$$ where
$$m_j = \begin{cases}   \frac1{\mu_j} \ \  {\text{ if}} \ \  \mu_j \neq 0\\
  \infty \ \ {\text{ if}} \ \   \mu_j = 0.
\end{cases} $$
\ed

 If none of the $m_j$ is
infinity, we say that $M$ is of {finite multitype at $p$}.

Clearly, since the definition of multitype  includes all
distinguished weights, the infimum is a { biholomorphic invariant}.

{Coordinates corresponding to the multitype weight $\Lambda_M$, in
which the local description of $M$ has form (\ref{1}), with $P$
being  $\Lambda_M$-homogeneous, are called
{ multitype coordinates}.}

If $M$ is of
finite multitype at $p$, the infimum in \eqref{2.6}  is
attained, which implies that  multitype coordinates do exist (\cite{C}, \cite{Ko1}).

\bd  Let $M$ be given by (2.3). 
We define a
 model hypersurface $M_P$ associated
to $M$ at $p=0$ by 
\begin{equation} M_P = \{(z,w) \in \mathbb C^{3}\ | \
 v  = P \zz \}. \label{22}\end{equation}\ed


%
%
%


%

Next let us  recall the following definitions.

\bd{Let $X$ be  a holomorphic vector field
of the form
\begin{equation}
X = \sum_{j=1}^2 f^j(z,w) \partial_{ z_j} + g(z,w)\partial_{w}.
\end{equation}
We say that  $X$  is rigid if $f^1,  f^2, g $ are all independent of the variable  $w$.}
\ed

Note that
the  rigid vector field  $\ W,$ of homogeneous weight
$-1,$   given by
\begin{equation}\label{gentil.1}
 W = \partial_{w}
\end{equation}
lies in $\aut (M_P,0).$ We will denote by  $E$  the weighted
homogeneous vector field of weight $0$ defined  by
\begin{equation}
 E = \sum_{j=1}^2 \mu_j z_j \partial_{ z_j} + w\partial_{w},
\end{equation}
i.e.  the weighted Euler field.
 Note that by definition of $\mu_j,$ E is a non rigid vector field  lying in $\aut (M_P,0).$

We can divide homogeneous rigid  vector fields into three types, and introduce
the following terminology.

\bd
Let $X \in \aut(M_P,0)$ be a rigid weighted homogeneous vector
field. $X$ is called
\begin{enumerate}
\item a { shift} if the weighted degree of $X$ is less than zero;
\item a { rotation }
if the weighted degree of $X$ is equal to  zero;
\item  a
 { generalized rotation}  if the weighted degree of $X$ is bigger than
 zero
\end{enumerate}
\ed

Notice that $X\in \aut(M_P, 0) $
is a generalized rotation if and only if it has positive weighted degree
and commutes with $W$. In other words, generalized rotations are precisely the elements of $\5g_c$. 

\section{Generalized rotations}

In this section we study nonlinear infinitesimal  CR automorphisms of hypersurfaces of finite multitype and  derive an
 explicit description of  all models which 
 admit a generalized rotation.

 \bl {Let $Y=f_1\frac{\partial }{\partial z_1} +f_2\frac{\partial }{\partial z_2} $ be a weighted 
homogeneous holomorphic  vector field of weighted degree $> 0$.  Then the space of weighted homogeneous polynomials  
in $z$ of a given weighted degree $\nu$
annihilated by $X$ has complex dimension at most one.}
\el

\bpf
%
%
%
%
First we claim that 
Y cannot be a multiple of the Euler field. Indeed, if $Y = h E$ with $h$ holomorpic and nonconstant, 
then Y(P)= 0  has no homogeneous  solution, since $\Re Y(P) = \Re h P \neq 0$.
Hence there exists
a point $q$ 
such that $Y(q)$ is not a multiple of the Euler field. 
This point lies on a uniquely determined complex curve 
$ z_1^{m_1} = c z_2^{m_2}$, and $Y$ is transverse to this curve in a neighbourhood of $q$. 
By homogeneity, on this curve $P (z_1, z_2)$ is determined up to a multiplicative complex constant.
Fixing this constant, by the uniqueness property for solutions of 
complex ODEs (\cite{IJ}), the equation $Y(P) = 0 $ determines $P$ uniquely in a neighbourhood of $q$, 
Since $P$ is a polynomial, if it exists, it is determined uniquely. 
Hence the space of solutions of $Y(P)=0 $ is at most 
one dimensional. 


%
%

\epf

\bl \label{do0} Let  $V_n, \ n\in \mathbb N$, be the  space
\begin{equation}
V_n=\{X |Y^n(X)=0                 \},
\end{equation}
where $X$ is a holomorphic polynomial of a given  constant weighted length  and $Y$ is a weighted homogeneous holomorphic vector field.  Then 
\begin{equation}
\dim V_n \le n.
\end{equation}
Moreover, when $d_n = \dim V_n >0,$  one can choose a basis for $V_n$ of the form
\begin{equation}\label{do} \{ {F_s^n}, \ s= {1},2, \dots ,{d_n}    
| \ \  Y^{d_n}({F_{d_n}^n}) =0, \ Y^{{d_n}-1}({F_{d_n}^n}) \ne 0, Y^{d_n-1}({F_{s}^n})=0,   s= {1},2, \dots ,{d_n-1}  \}
\end{equation}
\el
\bpf
We prove the lemma by induction. The case $n=1$ is a direct application of the previous Lemma. Suppose now that the lemma is true for $n$ and prove it for $n+1.$ We have 
\begin{equation}
V_{n+1}= \{X |Y^{n+1}(X)=0\}=\{X |Y^n( Y(X))=0\}.   
\end{equation}
By induction, we obtain that
\begin{equation}\label{do3}
Y(X)\in \span [ 
 {F_s^n}, \        |  \  Y^{d_n}({F_{d_n}^n}) =0, \\ \ Y^{{d_n}}({F_{d_n}^n}) \ne 0, Y^{d_n-1}({F_{s}^n})=0, s= {1}, \dots ,{d_n-1}   
]
\end{equation}
which implies that 
\begin{equation}
\dim V_{n+1}\le n+1.
\end{equation}
After performing  a linear combination of the solutions $X$ of \eqref{do3}, we may satisfy \eqref{do}.
\epf
\bt 
\label{2piva} Let 
 $M_P$  be  given by  \eqref{model} admitting a generalized rotation
 $Y.$  
Then $P$ can be decomposed in the following way

\begin{equation}P=\sum_{j=1}^M T_j,
\end{equation} where each $T_j$ is given by 
\begin{equation}\label{Barou2}
 T_j =\Re (\sum_{k=1}^{N_j}
  {U_j^{k}}
  {\overline {{V_j^{N_j -k +1}}}}), 
  \end{equation}
where $\{{ {{U_j^{1}}, \dots, {U_j^{N_j}} }}\}$  and  $\{{ {{V_j^{1}}, \dots, {V_j^{N_j}} }}\}$ are a symmetric pair of $Y-$ chains.

\et
\bpf
Let
\begin{equation}
P = \sum_{k=1}^l P_k,
\end{equation}
where $P_1 \neq 0, P_l \neq 0$, be the bihomogeneous expansion of $P$. Each $P_j$ is weighted homogeneous with respect 
to $z$ of weighted degree $c_j$  where $c_1< c_2 < \dots <c_l.$

We may write 
\begin{equation}
P_1 =\sum_{j=1}^{r} {S_j^{c_1}}{\overline {{S_j^{\hat c_1}}}},
\end{equation}
with $r$ minimal. Note that $c_1 + \hat c_1 =1.$ We claim that $r=1.$ Since $Y$ is a generalized rotation, we must have 
\begin{equation}\label{do1} \overline Y(\sum_{j=1}^{r} {S_j^{c_1}}{\overline {{S_j^{\hat c_1}}})=
\sum_{j=1}^{r} {S_j^{c_1}} \overline Y ({\overline {S_j^{\hat c_1}}}} )= 0.
\end{equation}
Since $r$ is minimal, this implies that
\begin{equation}\label{do2}
 Y ( S_j^{\hat c_1})= 0
 \end{equation}
 for all $j$. 
 Using Lemma~\ref{do0}, we conclude
 that 
\begin{equation} S_j^{\hat c_1} \in [S_1^{\hat c_1}  ]
\end{equation} for all  j.
  We may then write $P_1$ as 
\begin{equation}
P_1=Q_1^{c_1} \overline{Q_1^{\hat c_1}}.
\end{equation}
Hence, $r=1$ and the claim is proved.
We  write  now 
\begin{equation}
 P_k =\sum_{j=1}^{r_k} {S_j^{c_k}}{\overline{S_j^{\hat c_k}}}
\end{equation}
with $r_k$ minimal.

We claim that $P_k$ can be rewritten as 
\begin{equation}\label{dodo}
P_{k} =Q_k^{c_k} \overline{Q_k^{\hat c_k}}
+\tilde {P_{k}}
\end{equation}
so that there is a $d_k \leq k$ such that 
\begin{equation}
  \overline{ {Y}^{d_k}}(\overline{Q_k^{\hat c_k}})=0, \ \overline{ {Y}^{d_k-1}}
  (\overline{Q_k^{\hat c_k}})\ne 0,\  \overline{{Y}^{d_k-1}}(\tilde {P_{k} })=0.
\end{equation}

We prove the claim by induction. The case $k=1$ has just been proved.
%

Suppose  by induction that \eqref{dodo} holds for $k.$  Since $Y$ is a generalized rotation, we have
\begin{equation}\label{do5} Y({{Q}_k^{c_k}}){\overline
 {{{Q}_k^{\hat c_k}}} + Y(\tilde {P_{k} }) +\sum_{j=1}^{r_{k+1}} {S_{k+1}^{c_{k+1}}}
 \overline Y({\overline {S_{k+1}^{\hat c_{k+1}}}}})= 0.
\end{equation}
Applying ${\overline Y}^{d_k}$ to \eqref{do5}, we get
\begin{equation}\label{do6}
\sum_{j=1}^{r_{k+1}} S_{j}^{c_{k+1}} \overline{ Y^{d_{k}+1}}
(\overline{{S_{j}^{\hat c_{k+1}}}})= 0.
\end{equation}
Since $r_{k+1}$ is minimal,
\begin{equation}
 \overline{ Y^{d_{k}+1}}
(\overline{{S_{j}^{\hat c_{k+1}}}})=0
\end{equation}
for all $j$.
 Using Lemma~\ref{do0}, we obtain that $r_{k+1}\le d_k + 1 \leq k+1.$  
Using \eqref{do},  we may then  rewrite $P_{k+1}$  in the form given by \eqref{dodo}. The claim is then proved.
\\
Let $N_1\leq l $ be minimal such that 
$$   Y({{Q}_k^{c_k}})\ne 0,\   k=1, \dots, N_1-1,  \ \  Y({{Q}_{N_1}^{c_{N_1}}})= 
0.\}$$
We consider the following set $E_1$ given by
\begin{equation}\label{no}
E_1=\{{{Q}_k^{c_k}}{\overline {{{Q}_k^{\hat c_k}}}}, \ \ k=1,
 \dots,  N_1 \} 
\end{equation}
Note that this set is not empty since $ Y({Q_l^{c_l}})= 0$.

We claim that the following holds for every element of $E_1.$
\begin{enumerate}
\item $d_k=k, \ k=1, \dots, N_1. $   \\ 
\item $Y({{Q_k}^{c_k}})= a_{k} {{Q}_{k+1}^{c_{k+1}}} ,$\\
\item $ Y({{{Q_{k+1}^{\hat c_{k+1}}}}})=b_{{
{k+1}}} {
{{{Q_k}^{\hat c_k}}}} +R_k,$
where 
${Y}^{k-1}({R_{k} })=0.$\\
\end{enumerate}

We show that $d_k=k$ using  induction as above. Indeed, suppose that this is true for $k < N_1-1$ 
and show that it is also true for $k+1.$ Using the fact that $Y$ is a generalized rotation, we have as in \eqref{do5}

\begin{equation}\label{do7} Y({Q_k^{c_k}){\overline {{Q_k^{\hat c_k}}}} + Y(\tilde P_{k} ) +
({{Q}_{{k+1}}^{c_{k+1}}})\overline Y ({\overline {{Q}_{{k+1}}^{\hat c_{k+1}}}}}) + \overline Y(\tilde P_{k+1} )=0.
\end{equation}
Applying $ \overline Y^{k-1}$ to \eqref{do7}, we obtain 
\begin{equation}\label{do9} Y({Q_k^{c_k}})  \overline Y^{k-1}({\overline {{Q_k^{\hat c_k}}}}) 
+ ({{Q}_{{k+1}}^{c_{k+1}}})\overline Y^k ({\overline {{{Q}_{{k+1}}^{\hat c_{k+1}}}}}) =0.
\end{equation}
Hence, using \eqref{do9}, $d_{k+1}=k+1$ by definition of $E_1,$ and hence
\begin{equation}
Y({Q_k^{c_k}})= a_{k} {Q_{k+1}^{c_{k+1}}}.
\end{equation}
\begin{equation}
 Y^k({{{{Q}_{{k+1}}^{\hat c_{k+1}}}}})=b_{{
 {k+1}}}Y^{k-1}{{{{Q}_k^{\hat c_k}}}}, 
\end{equation}
which implies 
\begin{equation}\label{do10}
 Y^{k-1}(Y({{{{Q}_{k+1}^{\hat c_{k+1}}}}})-b_{{
 {k+1}}}{{{{Q}_k^{\hat c_k}}}})=0,
\end{equation}
and hence
\begin{equation}\label{do11}
 Y({{{{Q}_{k+1}^{\hat c_{k+1}}}}})=b_{{
 {k+1}}} { {{{Q}_k^{\hat c_k}}}} +R_k,
\end{equation}
where 
$ {Y}^{k-1}({R_{k} })=0.$ This achieves the proof of the claim.
Using \eqref{do11} and \eqref{dodo}, we may then assume without loss of generality that $R_k=0.$
We define the  chains by putting 
\begin{equation}\label{no1}
\begin{cases}
 {U}^{k}_{1}  :={Q}^{c_k}_{k},\\
V_1^{{{k}}}  :={Q}_k^{\hat c_{N_1-k+1}},\\

\end{cases}
\end{equation}
It follows from the above  properties of $E_1$ that $U_1^k$ and $V_1^k$
form a chain.   
In other words, we may write
\begin{equation}
P=\Re (\sum_{k=1}^{N_1}
  {U_1^{k}}
  {\overline {{V_1^{N_1 -k +1}}}}) + {\hat P}, \ \ \ k=1, \dots, N_1,
\end{equation}

It follows from \eqref{do9} that $Y$ is a generalized rotation for  $$ \Im w =\Re (\sum_{k=1}^{N_1}
  {U_1^{k}}
  {\overline {{V_1^{N_1 -k +1}}}}).$$
  It follows from \eqref{do9} that $a_k = - \bar b_{k+1}$, which means 
that the $U$ and $V$ are a pair of symmetric chains.
Hence $Y$ is a generalized rotation also for $\hat P$. We can repeat the above argument for 
$\hat P$ and in a finite number of steps we reach conclusion of the theorem.

\end{proof}

Note that symmetric chains and pairs of chains of any length can arise. 
\be
Let 
$$ Y = z_1^2 \d_{z_1} - z_1 z_2 \d_{z_2}. $$
Given three integers $1\leq l \leq m \leq n$
we first define
$$U^l = z_1^nz_2^n .$$
We can build a symmetric $Y$-chain by setting 
 $U^j = c_j z_1^{n-l+j}z_2^n $ for $j = 1, ..., l-1$
for suitable constants $c_j$.
Analogously, setting in addition  
$$V^l = z_1^mz_2^m $$
we can get in the same way  a pair of symmetric Y-chains of arbitrary length $l$.
\ee

\section{Proofs of the main results
}
In this section we complete the proofs of the results stated in the introduction.
The first part of Theorem \ref{pivo} has been already proved in Section 3 (as Theorem \ref{2piva}).
The second, converse part of the statement is immediate to verify.
%

In order to prove Theorem \ref{7piv}, we combine Theorem \ref{pivo} with Theorem 4.7 and Theorem 6.2 from \cite{KMZ}. 
They imply that on a smooth hypersurface of finite Catlin multitype 
2-jets are always sufficient to determine an element from $\5g_1$ and $\5g_n$. 

{\it Proof of Theorem \ref{3piva}}.
\ 
In the notation  of Theorem \ref{2piva}, we set
\begin{equation}
 K =  2 \sum_{j=1}^{M} N_j+1.
\end{equation}
We define a hyperquadric in  $\mathbb C^{K+1}$ by
\begin{equation}
 \Im \eta = \Re \sum_{j=1}^{M}\sum_{k=1}^{N_j}
   \zeta_{j,k}
  \overline{\zeta'_{j, N_j -k +1}}, 
\end{equation}
and consider the mappping  $\mathbb C^3 \to \mathbb C^{K+1}$
given by $\eta = w$ and 
\begin{equation}
 \zeta_{j,k} = U_j^{k}(z_1, z_2).
\end{equation}
and
\begin{equation}
 \zeta'_{j,k} = V_j^{k}(z_1, z_2).
\end{equation}
It is immediate to verify that the automorphism $Y$ of $M_P$ is $f$-related to the automorphism of this hyperquadric, defined by 
\begin{equation}
 Z =  \sum_{j=1}^{M}\sum_{k=2}^{N_j}
 c_{{k-1}, j} \;  \zeta_{j,{k}} \d_{ \zeta_{j,k-1 }} + d_{{k-1}, j} \;  {\zeta}'_{j,{k}} \d_{ {\zeta}'_{j,k-1 }}.
\end{equation}
Indeed, the condition for $f$-related vector fields becomes exactly the chain condition 
\eqref{yuv}-\eqref{dcj}.



%


\begin{thebibliography}{BER96b}

\itemsep=2pt

\bibitem{BEH}
Baouendi, M. S.
 Ebenfelt, P.,
Huang, X.
 {\it Super-rigidity for CR embeddings of real hypersurfaces into hyperquadrics}
Adv. Math.
{\bf 219}
(2008),
1427--1445.

\bibitem{BP}
{Bedford, E.,  Pinchuk, S. I.,
 {\it Convex domains with noncompact groups of automorphisms},
  {Mat. Sb.},
    \textbf{185}
     (1994), 3--26.

\bibitem{BES}
Beloshapka, V. K., Ezhov, V. V.,  Schmalz, G.,
{\it Holomorphic classification of four-dimensional surfaces in
              {$\mathbb C^3$}},
 Izv. Ross. Akad. Nauk Ser. Mat.},
 \textbf{72} (2008), 3--18.

\bibitem{BK}
Beloshapka, V. K., Kossovskiy, I. G.,
{\it Classification of homogeneous {CR}-manifolds in dimension 4},
J. Math. Anal. Appl.,
  \textbf{374}
 (2011),
 655--672.

\bibitem{BG}
Bloom, T.,  Graham, I., {\it On "type"  conditions for
generic real submanifolds of $  C\sp{n}$},  Invent. Math.  \textbf{40}
(1977),  217--243.

\bibitem{BG1}  Bloom, T., Graham, I. {\it A geometric characterization of 
points of type $m$ on real submanifolds of $C\sp{n}$.}
  J. Differential Geometry  {\bf 12} (1977), 171--182.



\bibitem{BFG} Beals, M., Fefferman, C., Graham R.,
{\it Strictly pseudoconvex domains in $\mathbb C^n$},
 Bull.\ Amer.\ Math.\ Soc.\ (N.S.)
 \textbf{8} (1983), 125--322.
 
 \bibitem{Ca}  Cartan, E., 
 \textit{Sur la g\'eom\'etrie pseudo-conforme des hypersurfaces de deux
variables complexes, I }, Ann. Math. Pura Appl.  \textbf{11}
(1932), 17--90.


\bibitem{C}
  Catlin, D.,
   {\it Boundary invariants of pseudoconvex domains},
Ann.\ Math.\ {\bf 120} (1984), 529--586.

\bibitem{C1}
  Catlin, D.,
   {\it Subelliptic estimates for $\bar \partial$-Neumann problem  on pseudoconvex domains},
Ann.\ Math.\ {\bf 126} (1987), 131--191.


\bibitem{CM} Chern, S.\ S.\, Moser, J., \textit{Real hypersurfaces in
complex manifolds},  Acta Math.\  \textbf{133} (1974),  219--271.

\bibitem{D}
 D'Angelo, J., {\em  Orders of contact, real
hypersurfaces and applications}, Ann.\  Math.\ {\bf  115} (1982),
615--637.

\bibitem{E} Ebenfelt, P., {\it New invariant tensors in CR structures and normal form for real hypersurfaces at a generic Levi degeneracy}, J. Diff. Geom. {\bf 50} (1998), 207--247.



\bibitem{EHZ}
Ebenfelt, P.,
Huang, X.,
Zaitsev, D.,
{\it The equivalence problem and rigidity for hypersurfaces embedded into
   hyperquadrics},
 Amer. J. Math.
{\bf 127}
(2005),
169--191.

\bibitem{ELZ} Ebenfelt, P., Lamel, B., Zaitsev, D., {\it Degenerate real hypersurfaces in $\mathbb C^2$ 
with few automorphisms}, Trans. Amer. Math. Soc. {\bf 361} (2009), 3241--3267.

\bibitem{ELZ2} Ebenfelt, P., Lamel, B., Zaitsev, D., {\it
Finite jet determination of local analytic CR automorphisms and
   their parametrization by 2-jets in the finite type case},
  {Geom. Funct. Anal.},
   {\bf 13},
   (2003),
   3, 546--573.



\bibitem{FK}
Fels, G., Kaup, W., {\it Classification of Levi degenerate homogeneous 
CR manifolds in dimension 5}, Acta Math. \textbf{201} (2008),  1--82.

\bibitem{HY} Huang, X.,  and Yin, W., 
{\it A Bishop surface with a vanishing Bishop invariant}, 
Invent. Math. {\bf 176} (2009), 461--520. 

\bibitem{IJ} Ilyashenko, Y.,   Yakovenko, S., {\it Lectures on Analytic Differential Equations,} 
Graduate Studies in Mathematics, vol. 86, Amer. Math. Soc., 2007.

\bibitem{IK}
Isaev, A. V., Krantz, S. G.,
 {\it Domains with non-compact automorphism group: a survey},
  {Adv. Math.},
\textbf{146},
 (1999),
   1--38.

\bibitem{J}
Jacobowitz, H., {\it An Introduction to CR structures}, Mathematical Surveys and Monographs, 32.
Amer. Math. Soc., Providence, RI, 1990.

\bibitem{K} Kohn, J.\ J., \textit{Boundary behaviour of
$\bar \partial$ on weakly pseudoconvex manifolds of dimension
two},
 J.\ Differential  Geom.\  \textbf{6} (1972),  523--542.

%

\bibitem{KMZ} Kol\'a\v r, M., Meylan, F., Zaitsev, D.,  \textit{Chern-Moser operators and 
polynomial models
in CR geometry}, Adv.  Math. {\bf 263} (2014), 321--356.



\bibitem{KM3} Kol\'a\v r, M., Meylan, F., \textit{
Infinitesimal CR automorphisms of hypersurfaces of finite type in
   $\mathbb C^2$},
  Arch. Math., {\bf 47},
   (2011), 
   367--375.
\bibitem{Ko1} Kol\'a\v r, M., {\em The Catlin multitype and biholomorphic equivalence
of models}, Int. Math. Res. Not. IMRN {\bf 18} (2010), 3530--3548.

\bibitem{Ko1a} Kol\'a\v r, M., {\it Normal forms for hypersurfaces of finite type in $ \mathbb C^2$}, Math. Res. Lett., \textbf{12} (2005),  523--542.

\bibitem{FM}
Kol\'a\v r, M., Meylan, F., {\it Chern-Moser operators and weighted jet determination problems},
Geometric analysis of several complex variables and related topics, 75--88, Contemp. Math. 550, 2011.

\bibitem{KMp} Kol\'a\v r, M., Meylan, F., \textit{
 Automorphism groups  of Levi degenerate hypersurfaces in
   $\mathbb C^3$}, preprint

\bibitem{KL}
Kruzhilin, N. G., Loboda, A. V.,
 {\it Linearization of local automorphisms of pseudoconvex surfaces},
 Dokl. Akad. Nauk SSSR,
\textbf{271} (1983), 280--282.


\bibitem{Po} Poincar\'e, H.,  \textit{Les fonctions analytiques de
deux variables et la repr\'esentation conforme, } Rend. Circ. Mat.
Palermo \textbf{23} (1907), 185--220.


\bibitem{S} Stanton, N.,
{\it Infinitesimal CR automorphisms of real hypersurfaces}, Amer. J. Math. {\bf 118} (1996),  209--233.


\bibitem{V} Vitushkin, A.G., \textit{Real analytic
hypersurfaces in complex manifolds}, Russ. Math. Surv. \textbf{40}
(1985),  1--35.

\bibitem{W} Webster, S.M., \textit{On
the Moser normal form at a non-umbilic point}, Math. Ann.
 \textbf{233} (1978),  97--102.






%
%
%
%
%
%
%
%
%
%
%
%
%
%
%
%
%
%
%
%
 \end{thebibliography}
\end{document}